\documentclass[a4paper,10pt]{article}
\usepackage[english]{babel}
\usepackage[dvips]{graphicx}
\usepackage{float}
\usepackage{amsfonts}
\usepackage{amssymb}
\usepackage{amsmath}
\usepackage{geometry}
\newcommand{\G}{\Gamma}
\newcommand{\Z}{\mathbb{Z}}
\newcommand{\eps}{\varepsilon}
\newcommand{\cO}{\cal O}
\newcommand{\cE}{\cal E}

\title{On $d$-graceful labelings}

\author{Anita Pasotti \thanks{Dipartimento di Matematica,
Facolt\`a di Ingegneria, Universit\`a degli Studi di Brescia,
Via Valotti, 9, I-25133 Brescia, Italy. email: anita.pasotti@ing.unibs.it}}

\date{}
\newtheorem{defi}{Definition}[section]
\newtheorem{prop}[defi]{Proposition}

\newtheorem{ex}[defi]{Example}
\newtheorem{thm}[defi]{Theorem}

\begin{document}
\maketitle

\begin{abstract}
In this paper we introduce a generalization of the
well known concept of a graceful labeling.
Given a graph $\G$ with $e=d\cdot m$ edges, we call $d$-graceful labeling of $\G$
an injective function from  $V(\G)$ to the set $\{0,1,2,\ldots, d(m+1)-1\}$
such that $\{|f(x)-f(y)|\ |\ [x,y]\in E(\G)\} =\{1,2,3,\ldots,d(m+1)-1\}
 - \{m+1,2(m+1),\ldots,(d-1)(m+1)\}.$ In the case of $d=1$ and of $d=e$ we find the classical notion
 of a graceful labeling and of an odd graceful labeling, respectively.
 Also, we call $d$-graceful $\alpha$-labeling
 of a bipartite graph $\G$ a $d$-graceful labeling of $\G$ with the property that
 its maximum value on one of the two bipartite sets does not reach  its minimum value
on the other one. We show that these new concepts allow to obtain certain cyclic
graph decompositions.
We investigate the existence of $d$-graceful $\alpha$-labelings
for several classes of bipartite graphs, completely solving the problem for paths and stars
and giving partial results about cycles of even length and ladders.
\end{abstract}


\section{Introduction}
We will denote by $K_v$ and $K_{m \times n}$ the
\emph{complete graph on $v$ vertices} and the
\emph{complete $m$-partite graph with parts of size $n$}, respectively.
For any graph $\G$ we write $V(\G)$ for the set of its vertices,
$E(\G)$ for the set of its edges, and $D(\G)$ for the set of its \emph{di-edges},
namely the set of all ordered pairs $(x,y)$ with $x$ and $y$ adjacent vertices of $\G$.
If $|E(\G)|=e$, we say that $\G$ has \emph{size} $e$.

Let $\G$ be a subgraph of a graph $K$. A $\G$-\emph{decomposition} of $K$ is a set
of graphs, called \emph{blocks}, isomorphic to $\G$ whose edges partition
$E(K)$. In the case that $K=K_v$ one also speaks of a \emph{$\G$-system of order $v$}.
An \emph{automorphism} of a $\G$-decomposition $\cal D$ of $K$ is a bijection
on $V(K)$ leaving $\cal D$ invariant. A $\G$-decomposition of $K$ is said to be \emph{cyclic}
if it admits an automorphism consisting of a single cycle of length $|V(K)|$.
In this case for giving the set $\cal B$ of blocks it is enough to
give a complete system of representatives for the orbits of
$\cal B$ under the cyclic group. The blocks of such a system are usually called
\emph{base blocks}.\\
For a survey on graph decompositions we refer to \cite{BE}.

The concept of a {\it graceful labeling} of $\Gamma$  introduced by A. Rosa \cite{R}
is quite related to the existence problem for cyclic $\Gamma$-systems.
A \emph{graceful labeling} of a graph $\G$ of size $e$ is an injective map
$f$ from $V(\Gamma)$ to the set of integers $\{0,1,2,\ldots,e\}$
such that
$$ \{|f(x)-f(y)|\ |\ [x,y]\in E(\G)\}=\{1,2,\ldots,e\}.$$
The graph $\Gamma$ is said to be {\it graceful} if there exists a graceful labeling of it.
In the case that $\Gamma$ is bipartite and $f$ has the additional property that
its maximum value on one of the two bipartite sets does not reach its minimum value on the other one, one says that
$f$ is an {\it $\alpha$-labeling}.\\
Every graceful labeling of a graph $\Gamma$ of size $e$ gives rise to a cyclic
$\Gamma$-system of order $2e+1$ but an $\alpha$-labeling
of $\Gamma$ gives much more; it gives in fact a cyclic
$\Gamma$-system of order $2et+1$ for every positive
integer $t$ (see \cite{R}) and other kinds of graph decompositions
such as, for instance, a {\it bicyclic} $\Gamma$-decomposition of
$K_{2\times e}$ (see \cite{EV}) and a cyclic
$\Gamma$-decomposition of
$K_{m\times e}$ for every odd integer $m$ coprime with $e$ (see \cite{BG}).\\
For a very rich survey on
graceful labelings we refer to \cite{G}.

Many variations of graceful labelings are known. In particular Gnana Jothi \cite{GJ}
defines an \emph{odd graceful labeling} of a graph $\G$ of size $e$ as an injective
function $f: V(\G) \rightarrow \{0,1,2,\ldots,2e-1\}$ such that
$$ \{|f(x)-f(y)|\ |\ [x,y]\in E(\G)\}=\{1,3,5,\ldots,2e-1\}.$$
If such a function exists $\G$ is said to be \emph{odd graceful}.
She proved that every graph with an $\alpha$-labeling is also odd graceful, while the converse, in general, is not true.
One of the applications of these labelings is that trees of size $e$, with a suitable
odd graceful labeling, can be used to generate cyclic decompositions
of the complete bipartite graph $K_{2 \times e}$.
For results on odd graceful graphs see \cite{Ba, E, SAME}.

In this paper we propose the following new definition which is,
at the same time, a generalization of the concepts of a graceful labeling
and of an odd graceful labeling.
\begin{defi}
Let $\G$ be a graph of size $e$ and let $d$ be a divisor of $e$, say
$e=d\cdot m$. A $d$-\emph{graceful labeling} of $\G$
is an injective function $f:V(\G) \rightarrow \{0,1,2,\ldots, d(m+1)-1\}$
such that
\begin{align*}
\{|f(x)-f(y)|\ |\ [x,y]\in E(\G)\} &=\{1,2,3,\ldots,d(m+1)-1\}\\
&\quad  - \{m+1,2(m+1),\ldots,(d-1)(m+1)\}.
\end{align*}
\end{defi}
If $\G$ admits a $d$-graceful labeling
we will say that $\G$ is a $d$-\emph{graceful graph}.
Note that  we find again the concepts of a graceful labeling and of an
 odd graceful labeling in the two extremal cases of $d=1$ and $d=e$, respectively.
By saying that $d$ is \emph{admissible} we will mean that it is a divisor of $e$ and
hence that it makes sense to investigate the existence of a $d$-graceful labeling of $\G$.\\
\\
It is well known that the complete graph $K_v$ is graceful
if and only if $v\leq4$, see \cite{SWG}.
At the moment we are not able to extend this result to $d$-graceful labelings of the
complete graph with $d>1$.
%It is not hard to see that if there exists a $d$-graceful labeling of $K_v$ then
%$d\leq \frac{v}{2}$ and that such condition is necessary but not sufficient.
%For instance one can easily
%check that $K_4$ is not a $2$-graceful graph.
%So $K_4$ is a graceful, but not a $2$-graceful graph.
%On the contrary, $K_5$ is not a graceful graph, but it is a
%$2$-graceful graph. A $2$-graceful labeling of $K_5$ can be obtained
%labeling its vertices with the integers $\{0,1,4,9,11\}$.
We recall that a $(v,d,k,\lambda)$-\emph{relative difference set} (RDS) in a group $G$
is a $k$-subset $S$ of $G$ such that the list of all differences $x-y$ with $(x,y)$ an ordered pair of
distinct elements of $S$ covers $G-H$ exactly $\lambda$ times,
and no element of $H$ at all with $H$ a suitable subgroup of $G$ of order $d$
(see \cite{Pott} for a survey on this topic). A RDS as above
is said to be cyclic if $G$ is such.
It is immediate to recognize that every $d$-graceful labeling of $K_v$ can be seen as
a cyclic $(v(v-1)+2d,2d,v,1)$-RDS
though the converse is not true in general.
For instance, the labeling of the vertices of $K_5$ with the elements of the set
$S=\{0,1,4,9,11\}$ is $2$-graceful and we see that $S$ is a $(20,4,5,1)$-RDS.
Now note that $S'=\{0,1,3,11,20\}$ is also a cyclic $(20,4,5,1)$-RDS, but it is evident
that the labeling of the vertices of $K_5$ with the elements of $S'$ is not $2$-graceful.
The very few known results about cyclic RDSs do not allow us to say
much more about the $d$-gracefulness of the complete graph with $d>1$.\\
\\
The $\alpha$-labelings can be generalized in a similar way.
\begin{defi}
A $d$-\emph{graceful} $\alpha$-\emph{labeling} of a bipartite graph $\G$ is a
$d$-graceful labeling of $\G$  having the property that
its maximum value on one of the two bipartite sets does not reach its minimum value on the other one.
\end{defi}
In the next section we will see how these new concepts
can be used to obtain some graph decompositions. Then
we investigate the existence of $d$-graceful $\alpha$-labelings for several
classes of bipartite graphs. In particular, we prove that paths and stars admit a
$d$-graceful $\alpha$-labeling for any admissible $d$. We also present partial results about
$d$-graceful labelings of  even cycles and ladders.

\section{$d$-graceful labelings and graph decompositions}
It is known that graceful labelings are related to difference families,
see  \cite{AB}.
The concept of a \emph{relative difference set} has been generalized in \cite{B98}
to that of a \emph{relative difference family}. A further generalization has been
introduced in \cite{BP}.
\begin{defi}
Let $\G$ be a graph and let $d$ be a divisor of $v$.
A $(v,d,\G,1)$-\emph{difference family (DF} for short$)$ is a collection $\cal F$
of injective maps from $V(\G)$ to $\Z_v$ such that the list
\[\Delta{\cal F}=\{f(x)-f(y)\ |\ f\in {\cal F};\ (x,y)\in D(\G)\}\]
covers $\Z_v-\frac{v}{d}\Z_v$ exactly once while it does not contain any element of $\frac{v}{d}\Z_v$,
where $\frac{v}{d}\Z_v$ denotes the subgroup of $\Z_v$ of order $d$.
\end{defi}
If a $(v,d,\G,1)$-DF consists of a single map $f$ one says that $f$
is a \emph{difference graph}.
In \cite{BP} it is proved that
\begin{thm}\label{thm:DF}
If there exists a $(v,d,\G,1)$-DF then there exists
a cyclic $\G$-decomposition of $K_{\frac{v}{d} \times d}$.
\end{thm}
The next easy proposition establishes the link between $d$-graceful  labelings and DFs.
\begin{prop}
A $d$-graceful labeling of a graph $\G$ of size $e$
determines
%a $(\Z_{d(m+1)},(m+1)\Z_{d(m+1)},\G,2)$-DF and
a $\left(2d(\frac{e}{d}+1),2d,\G,1\right)$-DF.
\end{prop}
Proof. Let $f:V(\G)\rightarrow \left\{0,1,\ldots, d\left(\frac{e}{d}+1\right)-1\right\}$ be a $d$-graceful labeling of $\G$.
It is obvious that
%$\phi: x \in V(\G)\rightarrow f(x)\in \Z_{d(m+1)}$ is
%a $(\Z_{d(m+1)},(m+1)\Z_{d(m+1)},\G,2)$-DF and that
$\phi: x \in V(\G)\rightarrow f(x)\in \Z_{2d\left(\frac{e}{d}+1\right)}$ is
a $\left(2d(\frac{e}{d}+1), 2d,\G,1\right)$-DF.
\hfill $\Box$\\
\\
By the above considerations we can state
\begin{prop}\label{prop:d}
If there exists a $d$-graceful labeling of a graph $\G$ of size $e$ then there exists
%a cyclic $\G$-decomposition of $^2K_{(m+1)\times d}$ and
a cyclic $\G$-decomposition of $K_{\left(\frac{e}{d}+1\right)\times 2d}$.
\end{prop}

\begin{figure}[!pht]
\begin{center}
\includegraphics[width=0.2\textwidth]{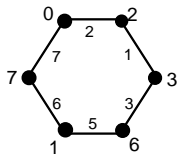}
\caption{A $2$-graceful labeling of $C_6$}
\label{2C6}
\end{center}
\end{figure}

\begin{ex}
Figure {\rm \ref{2C6}} shows a $2$-graceful labeling of $C_6$.
Hence $B=(0,2,3,6,1,7)$ is a base block of a cyclic $C_6$-decomposition of
$K_{4 \times 4}$.
\end{ex}
The following result shows that, as in the case of classical graceful labelings,
$d$-graceful $\alpha$-labelings are more powerful  than $d$-graceful
labelings.

\begin{thm}\label{prop:alpha}
If there exists a $d$-graceful $\alpha$-labeling of a graph $\G$ of size
$e$ then there exists
a cyclic $\G$-decomposition of $K_{\left(\frac{e}{d}+1\right)\times 2dn}$ for any integer $n\geq 1$.
\end{thm}
Proof. Let $f$ be a $d$-graceful $\alpha$-labeling of $\G$, so that we have
$\max f(X)<\min f(Y)$ where $X$ and $Y$ are the two parts of $\G$.
It is easy to see that
the set of maps $\{f_1,f_2,\ldots,f_n\}$
from $V(\G)$ to $\Z_{2dn\left(\frac{e}{d}+1\right)}$
defined by
$$f_i(x)=f(x)\quad \quad \forall x \in X,\ \forall i=1,\ldots, n$$
$$f_i(y)=f(y)+(i-1)(d+e)\quad \quad \forall y \in Y,\ \forall i=1,\ldots, n$$
is a $\left(2dn(\frac{e}{d}+1),2dn,\G,1\right)$-DF. Hence by Theorem \ref{thm:DF} there exists
a cyclic $\G$-decomposition of $K_{\left(\frac{e}{d}+1\right)\times 2dn}$ for any integer $n\geq 1$.
\hfill $\Box$

\begin{figure}[!pht]
\begin{center}
\includegraphics[width=0.2\textwidth]{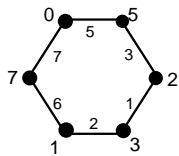}
\caption{A $2$-graceful $\alpha$-labeling of $C_6$}
\label{a2C6}
\end{center}
\end{figure}

\begin{ex}
Figure {\rm \ref{a2C6}} shows a $2$-graceful $\alpha$-labeling of
$C_6$ such that there exists a cyclic
 $C_6$-decomposition
of $K_{4\times 4n}$ for any integer $n$, $n \geq 1$.
The set of base blocks of such a decomposition is given by
$$ \{(0,5+8i,2,3+8i,1,7+8i)\ |\ i=0,\ldots, n-1 \}.$$
\end{ex}

\section{$d$-graceful $\alpha$-labelings of paths and stars}
In this section we will show that paths and stars have a
$d$-graceful $\alpha$-labeling for any admissible $d$. From now on, given two integers
$a$ and $b$, by $[a,b]$ we will denote the set of integers $x$ such that
$a\leq x\leq b$.\\
As usual we will denote by $P_{e+1}$ the path on $e+1$ vertices,
having size $e$. We recall that any path admits an $\alpha$-labeling, see \cite{R}.
\begin{thm}\label{path}
Given a positive integer $e$, the path $P_{e+1}$ has a $d$-graceful
$\alpha$-labeling for any admissible $d$.
\end{thm}
Proof. Let $e=d\cdot m$ and let $P_{dm+1}=(x_1 \sim x_2 \sim \ldots \sim  x_{dm+1})$.
We distinguish two cases depending on the parity of $m$.

\medskip
\noindent
Case 1: $m$ even.\\
Denote by ${\cO}=\{x_1,x_3,\ldots,x_{dm+1}\}$
and ${\cE}=\{x_2,x_4,\ldots,x_{dm}\}$
the bipartite sets of $P_{dm+1}$.
Let $f:V(P_{dm+1})\rightarrow [0,d(m+1)-1]$ be defined as follows:
$$f(x_{2i+1})=
i \quad \textrm{for}\ i \in \left[ 0,\frac{dm}{2} \right]$$
$$ f(x_{2i})=\left\{
\begin{array}{lll}
d(m+1)-i &\quad& \textrm{for}\ i \in \left[ 1,\frac{m}{2} \right]\\[2pt]
d(m+1)-i-1 &\quad& \textrm{for}\ i \in \left[\frac{m}{2}+1,m \right]\\[2pt]
d(m+1)-i-2 &\quad& \textrm{for}\ i \in \left[ m+1,\frac{3m}{2} \right]\\[2pt]
\vdots&\quad & \vdots\\
d(m+1)-i-(d-1) & \quad & \textrm{for}\ i \in \left[ \frac{(d-1)m}{2}+1,\frac{dm}{2} \right].
\end{array}
\right.$$
One can check that
\begin{align*}
 f({\cO})&=\left[0,\frac{dm}{2}\right]\\
f({\cE})&=\left[\frac{dm}{2}+1,\frac{(d+1)m}{2}\right]\cup\left[\frac{(d+1)m}{2}+2,\frac{(d+2)m}{2}+1\right]\cup\\
&\cup \left[\frac{(d+2)m}{2}+3,\frac{(d+3)m}{2}+2\right]\cup\ldots \cup\left[d(m+1)-\frac{m}{2},d(m+1)-1\right].
\end{align*}
Hence $f$ is injective and $\max f({\cO})<\min f({\cE})$. Now set
\begin{align}\label{epsrho}
\eps_i=|f(x_{2i-1})-f(x_{2i})|, \quad \rho_i=|f(x_{2i})-f(x_{2i+1})|\quad \textrm{for}
\ i=1,\ldots, \frac{dm}{2}.
\end{align}
By a direct calculation, one can see that
\begin{align*}
& \left\{\eps_i, \rho_i\ |\ i=1,\ldots, \frac{m}{2}\right\} =[(d-1)(m+1)+1, d(m+1)-1];\\
& \left\{\eps_i, \rho_i\ |\ i=\frac{m}{2}+1,\ldots, m\right\} =[(d-2)(m+1)+1, (d-1)(m+1)-1];\\
%& \left\{\eps_i, \rho_i\ |\ i=m+1,\ldots, \frac{3m}{2}\right\} =[e+d-3m-2,e+d-2m-3];\\
&\ldots\\
& \left\{\eps_i, \rho_i\ |\ i=\frac{(d-2)m}{2}+1,\ldots, \frac{(d-1)m}{2}\right\} =\left[m+2,2(m+1)-1\right];\\
& \left\{\eps_i, \rho_i\ |\ i=\frac{(d-1)m}{2}+1,\ldots, \frac{dm}{2}\right\} =\left[1,m\right].
\end{align*}
This concludes case 1.

\medskip
\noindent
Case 2: $m$ odd.\\
We divide this case in two subcases depending on the parity of $d$.\\
Case 2a): $d$ odd.
In this case the bipartite sets are ${\cO}=\{x_1,x_3,\ldots,x_{dm}\}$
and ${\cE}=\{x_2,x_4,\ldots,x_{dm+1}\}$.
Let $f:V(P_{dm+1})\rightarrow [0,d(m+1)-1]$ be defined as follows:
$$ f(x_{2i+1})=\left\{
\begin{array}{lll}
i &\quad& \textrm{for}\ i \in \left[ 0,\frac{m-1}{2} \right]\\[2pt]
i+1 &\quad& \textrm{for}\ i \in \left[\frac{m+1}{2},\frac{3m-1}{2} \right]\\[2pt]
i+2 &\quad& \textrm{for}\ i \in \left[ \frac{3m+1}{2},\frac{5m-1}{2} \right]\\[2pt]
\vdots & \quad & \vdots\\
i+\frac{d-1}{2} & \quad & \textrm{for}\ i \in \left[ \frac{(d-2)m+1}{2},\frac{dm-1}{2} \right].
\end{array}
\right.$$
$$ f(x_{2i})=\left\{
\begin{array}{lll}
d(m+1)-i &\quad& \textrm{for}\ i \in \left[ 1,m \right]\\[2pt]
d(m+1)-i-1 &\quad& \textrm{for}\ i \in \left[m+1,2m \right]\\[2pt]
d(m+1)-i-2 &\quad& \textrm{for}\ i \in \left[ 2m+1,3m \right]\\[2pt]
\vdots & \quad & \vdots\\
d(m+1)-i-\frac{d-1}{2} & \quad & \textrm{for}\ i \in \left[ \frac{(d-1)m}{2}+1,\frac{dm+1}{2} \right].
\end{array}
\right.$$
In this case it results that
\begin{align*}
 f({\cO})&=\left[0,\frac{m-1}{2}\right]\cup\left[\frac{m+3}{2},\frac{3m+1}{2}\right]\cup\\
 &\cup \left[\frac{3m+5}{2},\frac{5m+3}{2}\right]\cup\ldots
 \cup \left[\frac{(d-2)m+d}{2},\frac{d(m+1)}{2}-1\right]\\
f({\cE})&=\left[\frac{d(m+1)}{2},\frac{(d+1)(m+1)}{2}-1\right] \cup\\ & \cup \left[\frac{(d+1)(m+1)}{2}+1,\frac{(d+3)(m+1)}{2}-1\right]\cup\\
&\cup \left[\frac{(d+3)(m+1)}{2}+1,\frac{(d+5)(m+1)}{2}-1\right]\cup\ldots\\
& \ldots \cup\left[d(m+1)-m,d(m+1)-1\right].
\end{align*}
So, also in this case, $f$ is injective and $\max f({\cO})< \min f({\cE})$.
Note that, now, since $dm$ is odd
we have to consider
\begin{align*}
& \eps_i=|f(x_{2i-1})-f(x_{2i})|,\quad\quad \rho_i=|f(x_{2i})-f(x_{2i+1})|\quad \textrm{for}
\ i=1,\ldots, \frac{dm-1}{2},\\
& \eps_{\frac{dm+1}{2}}=|f(x_{dm})-f(x_{dm+1})|.
\end{align*}
By a long and tedious calculation, one can see that
\begin{align*}
 \left\{\eps_i, \rho_i\ |\ i=1,\ldots, \frac{m-1}{2}\right\} \cup \left\{\eps_{\frac{m+1}{2}}\right\} =&[(d-1)(m+1)+1,\\
 &d(m+1)-1];\\
\left\{\rho_{\frac{m+1}{2}}\right\} \cup \left\{\eps_i, \rho_i\ |\ i=\frac{m+3}{2},\ldots, m\right\}=&[(d-2)(m+1)+1,\\
&(d-1)(m+1)-1];\\
\ldots&
\end{align*}
%%& \left\{\eps_i, \rho_i\ |\ i=m+1,\ldots, \frac{3m}{2}\right\} =[e+d-3m-2,e+d-2m-3];\\
\begin{align*}
\left\{\rho_{\frac{(d-2)m+1}{2}}\right\} \cup \left\{\eps_i, \rho_i\ |\ i=\frac{(d-2)m+1}{2}+1,\ldots, \frac{(d-1)m}{2}\right\} =&\\
=[m+2, 2(m+1)-1]&;\\
\left\{\eps_i, \rho_i\ |\ i=\frac{(d-1)m}{2}+1,\ldots, \frac{dm-1}{2}\right\} \cup
\left\{\eps_{\frac{dm+1}{2}} \right\}=\left[1,m\right]&.
\end{align*}
This completes the proof of Case 2a).\\
Case 2b): $d$ even.
%Write ${\cO}=\{x_1,x_3,\ldots,x_{e}\}$
%and ${\cE}=\{x_2,x_4,\ldots,x_{e+1}\}$.
Let $f:V(P_{dm+1})\rightarrow [0,d(m+1)-1]$ be defined as follows:
$$ f(x_{2i+1})=\left\{
\begin{array}{lll}
i &\quad& \textrm{for}\ i \in \left[ 0,\frac{m-1}{2} \right]\\[2pt]
i+1 &\quad& \textrm{for}\ i \in \left[\frac{m+1}{2},\frac{3m-1}{2} \right]\\[2pt]
i+2 &\quad& \textrm{for}\ i \in \left[ \frac{3m+1}{2},\frac{5m-1}{2} \right]\\[2pt]
\vdots & \quad & \vdots\\
i+\frac{d}{2} & \quad & \textrm{for}\ i \in \left[ \frac{(d-1)m+1}{2},\frac{dm}{2} \right]
\end{array}
\right.$$
$$ f(x_{2i})=\left\{
\begin{array}{lll}
d(m+1)-i &\quad& \textrm{for}\ i \in \left[ 1,m \right]\\[2pt]
d(m+1)-i-1 &\quad& \textrm{for}\ i \in \left[m+1,2m \right]\\[2pt]
d(m+1)-i-2 &\quad& \textrm{for}\ i \in \left[ 2m+1,3m \right]\\[2pt]
\vdots & \quad & \vdots\\
d(m+1)-i-\frac{d-2}{2} & \quad & \textrm{for}\ i \in \left[ \frac{(d-2)m}{2}+1,\frac{dm}{2} \right].
\end{array}
\right.$$
Also here, reasoning as in Case 2a), it is easy to see that $f$ is a
$d$-graceful $\alpha$-labeling of $P_{dm+1}$.
\hfill $\Box$\\

\begin{figure}[!pht]
\begin{center}
\includegraphics[width=0.9\textwidth]{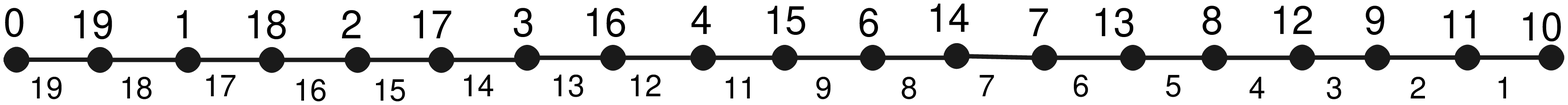}\\
\includegraphics[width=0.9\textwidth]{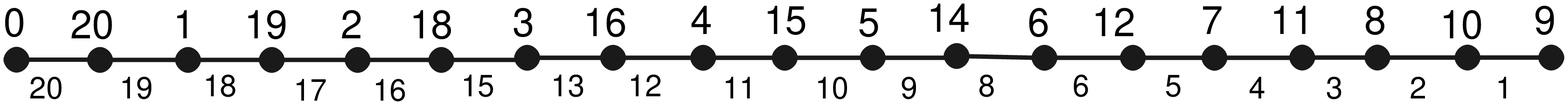}\\
\includegraphics[width=0.9\textwidth]{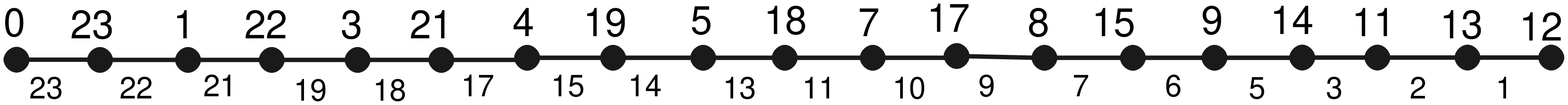}\\
\includegraphics[width=0.9\textwidth]{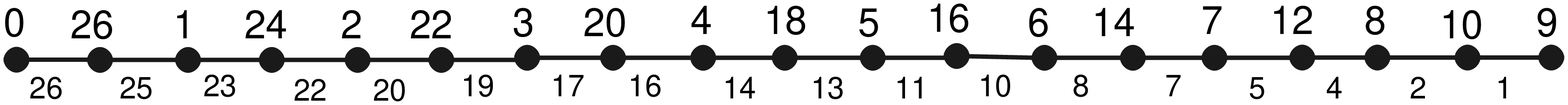}\\
\includegraphics[width=0.9\textwidth]{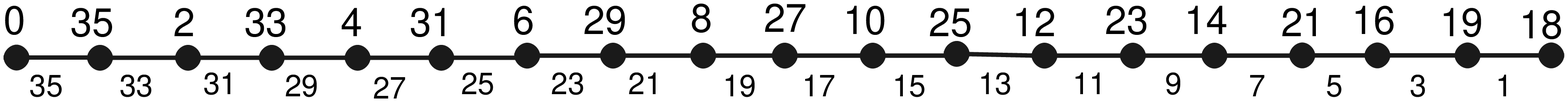}
\caption{$P_{19}$ }
\label{P19}
\end{center}
\end{figure}

\begin{prop}
There exists a cyclic $P_{e+1}$-decomposition of $K_{\left(\frac{e}{d}+1\right)\times 2dn}$
for any integers $e,n \geq 1$ and any divisor $d$ of $e$.
\end{prop}
Proof. It is a direct consequence of theorems \ref{prop:alpha} and \ref{path}.
\hfill$\Box$\\
\\
In Figure \ref{P19} we have the $d$-graceful $\alpha$-labelings of $P_{19}$
for $d=2,3,6,9, 18$, obtained following the construction
given in the proof of Theorem \ref{path}. Note that these are all the possible
$d$-graceful $\alpha$-labelings of $P_{19}$, apart from the already known case
of a classical $\alpha$-labeling $(d=1)$, \cite{R}.
\\
\\
We recall that the \emph{star} on $e+1$ vertices $S_{e+1}$ is the complete
bipartite graph with one part having a single vertex, called the
\emph{center} of the star, and the other part having $e$ vertices,
called \emph{external vertices}. It is obvious that any star has an $\alpha$-labeling.

\begin{thm}
Given a positive integer $e$, the star $S_{e+1}$ has a $d$-graceful
$\alpha$-labeling for any admissible $d$.
\end{thm}
Proof. Let $e=d\cdot m$. It is obvious that if we label the center
of the star with $0$ and the external vertices, arbitrarily, with the elements
of $\{1,2,3,\ldots,d(m+1)-1\} - \{m+1,2(m+1),\ldots,(d-1)(m+1)\}$, we have a
$d$-graceful $\alpha$-labeling of $S_{e+1}$.
\hfill $\Box$\\
\\
Applying Theorem \ref{prop:alpha} we have:
\begin{prop}
There exists a cyclic $S_{e+1}$-decomposition of $K_{\left(\frac{e}{d}+1\right)\times 2dn}$
for any integers $e,n \geq 1$ and any divisor $d$ of $e$.
\end{prop}

\section{$d$-graceful $\alpha$-labelings of cycles and ladders}
As usual, we will denote by $C_k$ the cycle on $k$ vertices. It is obvious that
$C_k$ is a graph of size $k$ and that it is bipartite if and only if
$k$ is even.\\
In \cite{R} Rosa proved that $C_k$ has an $\alpha$-labeling if and only if
$k\equiv$0(mod 4).

\begin{thm}\label{cycled2}
The cycle $C_{4k}$ has a $2$-graceful $\alpha$-labeling for every $k\geq1$.
\end{thm}
Proof. Let $C_{4k}=(x_1,x_2,\ldots,x_{4k})$ and
denote by ${\cO}=\{x_1,x_3,\ldots, x_{4k-1}\}$ and ${\cE}=\{x_2,x_4,\ldots,x_{4k}\}$
its bipartite sets of vertices.\\
Consider the map: $f: V(C_{4k}) \rightarrow \{0,1,\ldots, 4k+1\}$
defined as follows:
$$f(x_{2i+1})= i\quad  \textrm{for}\ i \in [0,2k-1]$$
$$ f(x_{2i})=\left\{
\begin{array}{lll}
4k+2-i &\quad& \textrm{for}\ i \in \left[ 1,k \right]\\
4k-i & \quad & \textrm{for}\ i \in \left[ k+1,2k \right].
\end{array}
\right.$$
We have
$$f({\cO})=[0,2k-1] \quad f({\cE})=\left[2k,3k-1\right]\cup\left[3k+2,4k+1\right].$$
In this way we see that $f$ is injective. Also, note that
$\max f({\cO})< \min f({\cE})$. Now set $\eps_i, \rho_i$ (for $i=1,\ldots,2k$) as in (\ref{epsrho}), where the indices are understood modulo $4k$.
One can easily check that
$$\left\{\eps_i, \rho_i\ |\ i=1,\ldots, k\right\}=[2k+2, 4k+1]$$
and
$$\left\{\eps_i, \rho_i\ |\ i=k+1,\ldots, 2k\right\}=[1,2k]. $$
Hence $f$ is a $2$-graceful $\alpha$-labeling of $C_{4k}$.
\hfill$\Box$

\begin{thm}\label{cycled4}
The cycle $C_{4k}$ admits a  $4$-graceful $\alpha$-labeling for every $k\geq1$.
\end{thm}
Proof.
Let  ${\cO}$ and ${\cE}$ be defined as in the proof
of the previous theorem. We are able to prove the existence of a $4$-graceful $\alpha$-labeling
of $C_{4k}$ by means of two direct constructions where we distinguish the two cases: $k$ even
and $k$ odd.

\medskip
\noindent
Case 1: $k$ even.\\
Consider the map: $f: V(C_{4k}) \rightarrow \{0,1,\ldots, 4k+3\}$
defined as follows:
$$f(x_{2i+1})= \left\{
\begin{array}{lll}
i &\quad& \textrm{for}\ i \in \left[ 0,\frac{3k}{2}-1 \right]\\[2pt]
i+1 & \quad & \textrm{for}\ i \in \left[ \frac{3k}{2},2k-1 \right]
\end{array}
\right.$$
$$ f(x_{2i})=\left\{
\begin{array}{lll}
4k+4-i &\quad& \textrm{for}\ i \in \left[ 1,\frac{k}{2} \right]\\[2pt]
4k+3-i & \quad & \textrm{for}\ i \in \left[ \frac{k}{2}+1,k \right]\\[2pt]
4k+1-i & \quad & \textrm{for}\ i \in \left[ k+1,2k \right].
\end{array}
\right.$$
It is easy to see that
\begin{align*}
f({\cO})& =\left[0,\frac{3k}{2}-1\right]\cup \left[\frac{3k}{2}+1,2k\right]\\
f({\cE})& =\left[2k+1,3k\right]\cup\left[3k+3,\frac{7k}{2}+2\right]\cup\left[\frac{7k}{2}+4,4k+3\right].
\end{align*}
Hence $f$ is injective and $\max f({\cO})<\min f({\cE})$.
Set $\eps_i$ and $\rho_i$ (for $i=1,\ldots,2k$) as in (\ref{epsrho}), where the indices are understood modulo $4k$.
It is not hard to see that
\begin{align*}
& \left\{\eps_i, \rho_i\ |\ i=1,\ldots, \frac{k}{2}\right\} =[3k+4, 4k+3];\\
& \left\{\eps_i, \rho_i\ |\ i=\frac{k}{2}+1,\ldots, k\right\} =[2k+3, 3k+2];\\
& \left\{\eps_i, \rho_i\ |\ i=k+1,\ldots, \frac{3k}{2}-1\right\} =[k+3, 2k];\\
& \left\{\eps_i, \rho_i\ |\ i=\frac{3k}{2},\ldots, 2k\right\} =[1,k]\cup \{k+2\}\cup \{2k+1\}.
\end{align*}
The thesis follows.

\begin{figure}[ht]
\begin{center}
\includegraphics[width=0.45\textwidth]{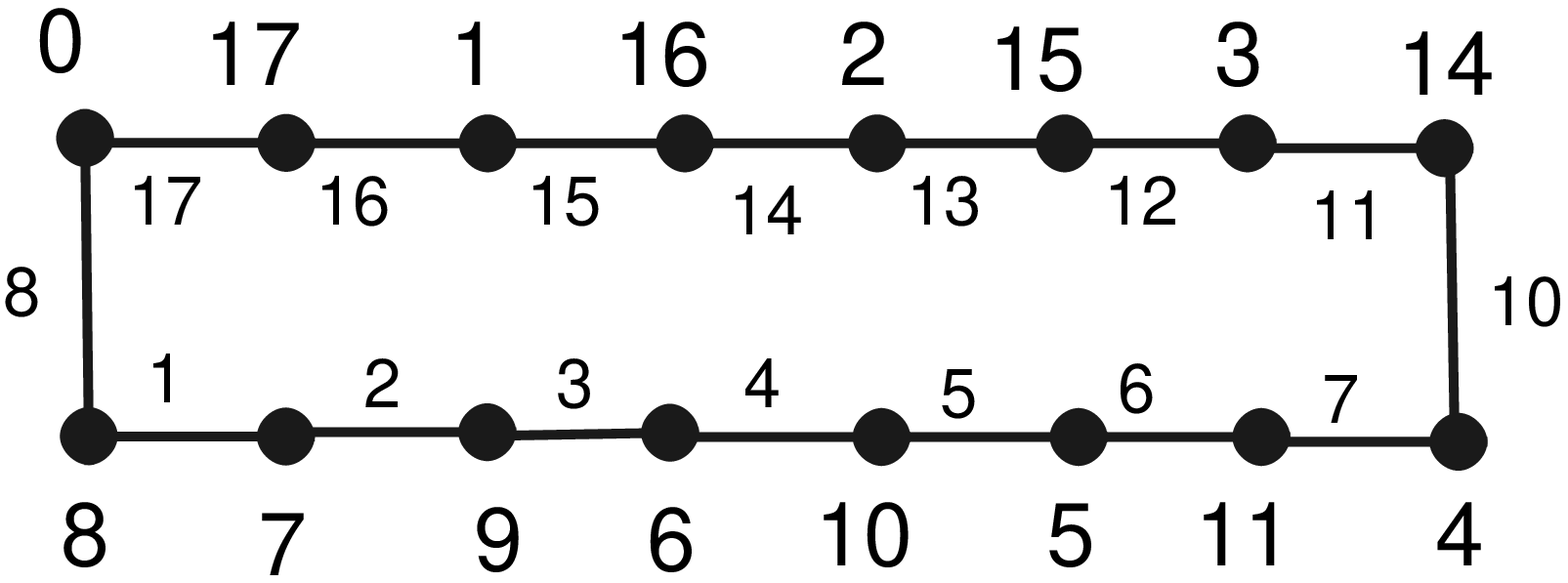}\quad
\includegraphics[width=0.45\textwidth]{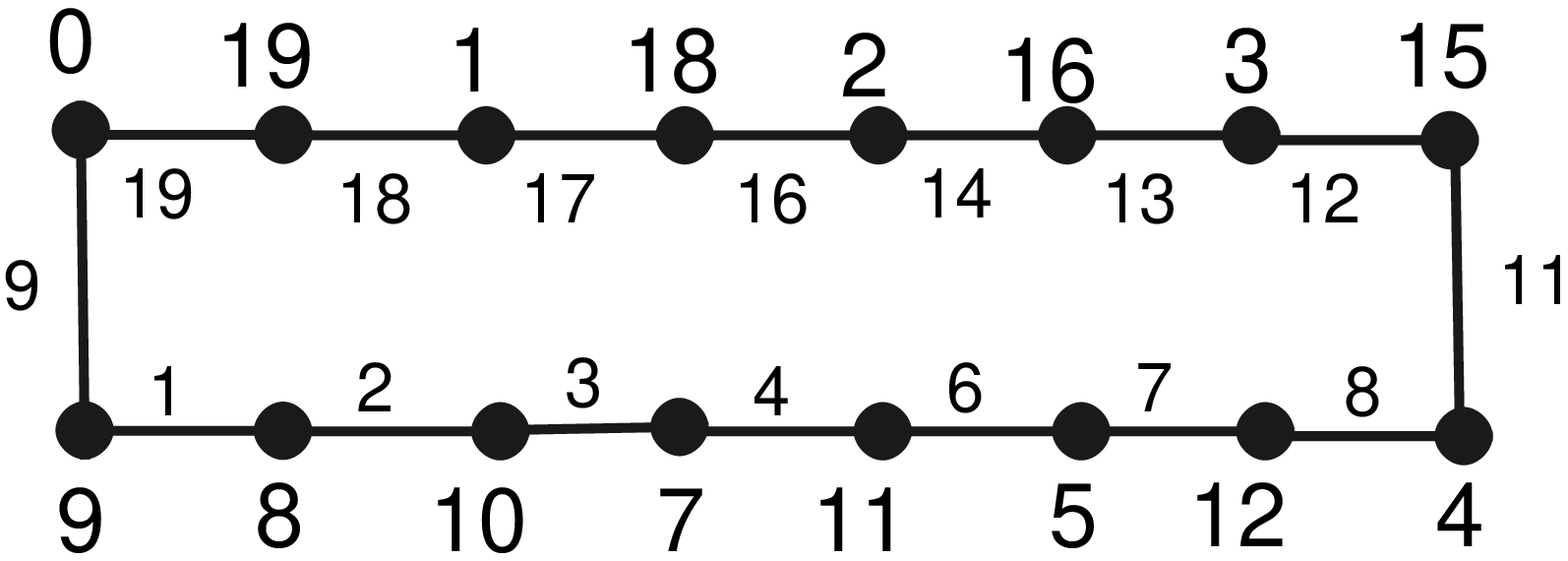}
\caption{$C_{16}$ }
\label{C16}
\end{center}
\end{figure}

\medskip
\noindent
Case 2: $k$ odd.\\
Let now $f: V(C_{4k}) \rightarrow \{0,1,\ldots, 4k+3\}$
be defined as follows:
$$f(x_{2i+1})= \left\{
\begin{array}{lll}
i &\quad& \textrm{for}\ i \in \left[ 0,\frac{k-1}{2} \right]\\[2pt]
i+1 & \quad & \textrm{for}\ i \in \left[ \frac{k+1}{2},2k-1 \right]
\end{array}
\right.$$
$$ f(x_{2i})=\left\{
\begin{array}{lll}
4k+4-i &\quad& \textrm{for}\ i \in \left[ 1,k\right]\\[2pt]
4k+2-i & \quad & \textrm{for}\ i \in \left[k+1,\frac{3k-1}{2} \right]\\[2pt]
4k+1-i & \quad & \textrm{for}\ i \in \left[ \frac{3k+1}{2},2k \right].
\end{array}
\right.$$
One can check that
\begin{align*}
f({\cO})& =\left[0,\frac{k-1}{2}\right]\cup \left[\frac{k+3}{2},2k\right]\\
f({\cE})& =\left[2k+1,\frac{5k+1}{2}\right]\cup\left[\frac{5k+5}{2},3k+1\right]\cup\left[3k+4,4k+3\right].
\end{align*}
Also in this case $f$ is injective and $\max f({\cO})<\min f({\cE})$.
By a direct calculation, one can see that
\begin{align*}
& \left\{\eps_i, \rho_i\ |\ i=1,\ldots, \frac{k-1}{2}\right\} =[3k+5, 4k+3];\\
& \left\{\eps_i, \rho_i\ |\ i=\frac{k+1}{2},\ldots, k\right\} =[2k+3, 3k+2]\cup\{3k+4\};\\
& \left\{\eps_i, \rho_i\ |\ i=k+1,\ldots, \frac{3k-1}{2}\right\} =[k+2,2k];\\
& \left\{\eps_i, \rho_i\ |\ i=\frac{3k+1}{2},\ldots, 2k\right\} =[1,k]\cup\{2k+1\}.
\end{align*}
The assertion then follows.
\hfill$\Box$\\
\\
The two figures in Figure \ref{C16} show the $2$-graceful $\alpha$-labeling and the $4$-graceful $\alpha$-labeling
of $C_{16}$ provided by the previous theorems.

\noindent
By virtue of Theorems \ref{prop:alpha}, \ref{cycled2}
and \ref{cycled4}, we have
\begin{prop}
There exists a cyclic $C_{4k}$-decomposition of $K_{(2k+1) \times 4n}$,
and a cyclic $C_{4k}$-decomposition of $K_{(k+1) \times 8n}$ for any integers $k,n\geq1$.
\end{prop}

\begin{thm}\label{C2k}
For any odd integer $k>1$, the cycle $C_{2k}$ is $2$-graceful.
\end{thm}
Proof. For $k=3,5,7$ we show directly a $2$-graceful labeling:
\begin{align*}
C_6&=(0,5,2,3,1,7),\\
C_{10}&=(0,11,1,3,7,4,5,10,2,9),\\
C_{14}&=(0,15,1,14,11,4,10,5,7,6,2,13,3,12).
\end{align*}
Let now $k$ be greater than $7$ and set $C_{2k}=(x_1,x_2,x_3,\ldots,x_{2k})$.
Let ${\cal O}=\{x_1,x_3,\ldots,x_{2k-1}\}$ and
${\cal E}=\{x_2,x_4,\ldots,x_{2k}\}$.
We will construct a $2$-graceful labeling of
$C_{2k}$ where we distinguish two cases according to whether $k\equiv 1$(mod 4) or
$k\equiv 3$(mod 4).

\medskip
\noindent
Case 1: $k\equiv 1$(mod 4).\\
In this case $k=2t+1$ with $t\geq 4$ even. Consider the map
$f: V(C_{2k})\rightarrow \{0,1,\ldots,2k+1\}$ defined as follows:
$$f(x_{2i+1})= \left\{
\begin{array}{lll}
i &\quad& \textrm{for}\ i \in \left[ 0,\frac{t}{2} \right]\\[2pt]
\frac{7}{2}t+3-i & \quad & \textrm{for}\ i \in \left[ \frac{t}{2}+1,t\right]\\[2pt]
3t+2-i & \quad & \textrm{for}\ i \in \left[ t+1,t+\lfloor\frac{t}{4}\rfloor\right]\\[2pt]
t+2+i & \quad & \textrm{for}\ i \in \left[ t+\lfloor\frac{t}{4}\rfloor+1,\frac{3}{2}t\right]\\[2pt]
i-t & \quad & \textrm{for}\ i \in \left[ \frac{3}{2}t+1,2t\right]
\end{array}
\right.$$
$$f(x_{2i})= \left\{
\begin{array}{lll}
4t+4-i &\quad& \textrm{for}\ i \in \left[1,\frac{t}{2} \right]\\[2pt]
\frac{1}{2}t+i & \quad & \textrm{for}\ i \in \left[ \frac{t}{2}+1,t+1\right]\\[2pt]
t+1+i & \quad & \textrm{for}\ i \in \left[ t+2,t+\lfloor\frac{t+2}{4}\rfloor\right]\\[2pt]
3t+2-i & \quad & \textrm{for}\ i \in \left[ t+\lfloor\frac{t+2}{4}\rfloor+1, \frac{3}{2}t  \right]\\[2pt]
5t+4-i & \quad & \textrm{for}\ i \in \left[ \frac{3}{2}t+1,2t+1\right].
\end{array}
\right.$$
If $t=4$ skip the third assignment of $f(x_{2i})$.\\
One can directly check that
$$f({\cal O})=\left[0,\frac{t}{2}\right]\cup \left[\frac{5}{2}t+3,3t+2\right] \cup
\left[2t+2-\left\lfloor\frac{t}{4}\right\rfloor, 2t+1\right]\cup$$
$$\cup \left[2t+3+\left\lfloor \frac{t}{4}\right\rfloor,\frac{5}{2}t+2 \right]
\cup \left[ \frac{t}{2}+1,t \right]$$

$$f({\cal E})=\left[\frac{7}{2}t+4,4t+3\right] \cup \left[t+1,\frac{3}{2}t+1\right] \cup
\left[2t+3,2t+1+\left\lfloor \frac{t+2}{4}\right\rfloor \right]\cup$$
$$\cup\left[\frac{3}{2}t+2,2t+1-\left\lfloor\frac{t+2}{4}\right\rfloor\right]\cup \left[ 3t+3,\frac{7}{2}t+3\right].$$
So $f$ is an injective function. Now set
$\eps_i$ and $\rho_i$ (for $i=1,\ldots, k$)
as in (\ref{epsrho})
where the indices are understood modulo $2k$.\\
By a direct calculation one can see that
\begin{align*}
& \left\{\eps_i,\rho_i\ |\ i=1,\ldots, \frac{t}{2}\right\}=[3t+4,4t+3]\\
& \eps_{\frac{t}{2}+1}=\frac{t}{2}+1\\
& \left\{\rho_{\frac{t}{2}+1}\right\}\cup \left\{\eps_i,\rho_i\ |\ i=\frac{t}{2}+2,\ldots, t\right\}\cup
\{\eps_{t+1}\}=[t+2,2t+1]\\
& \rho_{t+1}=\frac{t}{2}.
\end{align*}
We write the remaining differences and distinguish two subcases depending on the parity of $\frac{t}{2}$.
If $\frac{t}{2}$ is even, it results:
\begin{align*}
& \left\{\eps_i,\rho_i\ |\ i=t+2,\ldots,\frac{5}{4}t\right\}=\left[2,\frac{t}{2}-1\right]\\
& \eps_{\frac{5}{4}t+1}=1\\
& \left\{\rho_{\frac{5}{4}t+1}\right\} \cup \left\{\eps_i,\rho_i\ |\ i=\frac{5}{4}t+2,\ldots,\frac{3}{2}t\right\}
\cup \{\eps_{\frac{3}{2}t+1}\}=\left[\frac{t}{2}+2,t+1\right]\\
&\left\{\rho_{\frac{3}{2}t+1}\right\} \cup \left\{\eps_i,\rho_i\ |\ i=\frac{3}{2}t+2,\ldots,2t\right\}
\cup \left\{\eps_{2t+1}\right\}=[2t+3,3t+2]\\
& \rho_{2t+1}=3t+3.
\end{align*}
If $\frac{t}{2}$ is odd, we have:
\begin{align*}
& \left\{\eps_i,\rho_i\ |\ i=t+2,\ldots,\left\lfloor\frac{5}{4}t\right\rfloor\right\}\cup
\left\{\eps_{\left\lfloor\frac{5}{4}t\right\rfloor+1}\right\}=\left[2,\frac{t}{2}-1\right]\\
& \rho_{\left\lfloor\frac{5}{4}t\right\rfloor+1}=1\\
& \left\{\eps_i,\rho_i\ |\ i=\left\lfloor\frac{5}{4}t\right\rfloor+2,\ldots,\frac{3}{2}t\right\}
\cup \{\eps_{\frac{3}{2}t+1}\}=\left[\frac{t}{2}+2,t+1\right]\\
&\left\{\rho_{\frac{3}{2}t+1}\right\} \cup \left\{\eps_i,\rho_i\ |\ i=\frac{3}{2}t+2,\ldots,2t\right\}
\cup \left\{\eps_{2t+1}\right\}=[2t+3,3t+2]\\
& \rho_{2t+1}=3t+3.
\end{align*}
This concludes the proof of case 1.

\begin{figure}[!pht]
\begin{center}
\includegraphics[width=0.73\textwidth]{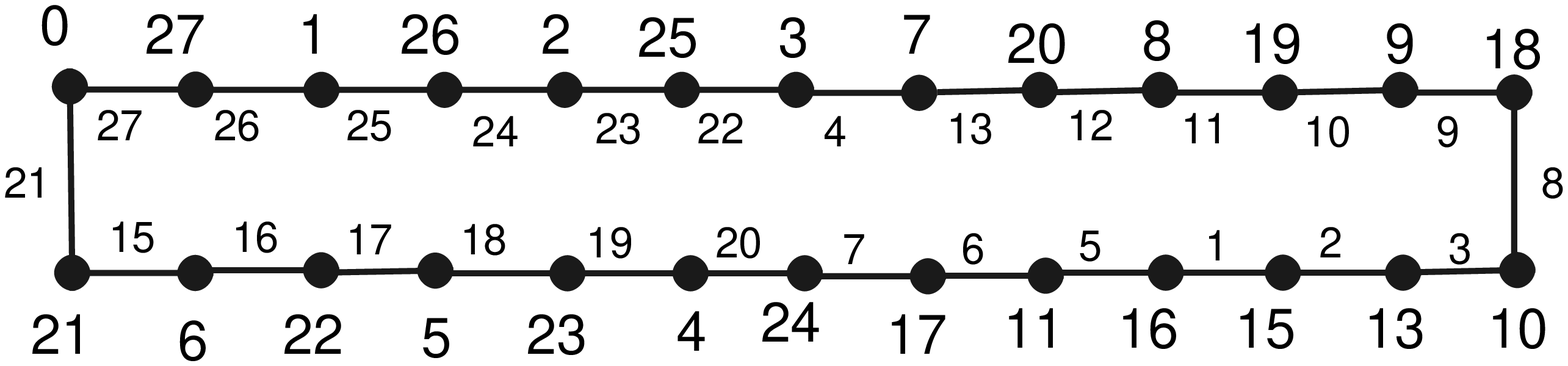}\\
\includegraphics[width=0.8\textwidth]{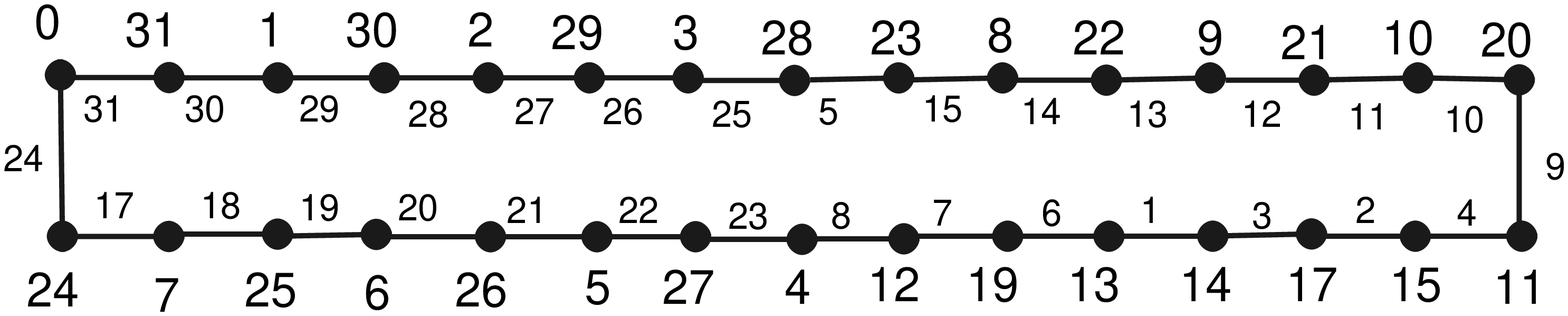}
\caption{$C_{26}$ and $C_{30}$}
\label{C18-22}
\end{center}
\end{figure}

\begin{figure}[h!]
\begin{center}
\includegraphics[width=0.45\textwidth]{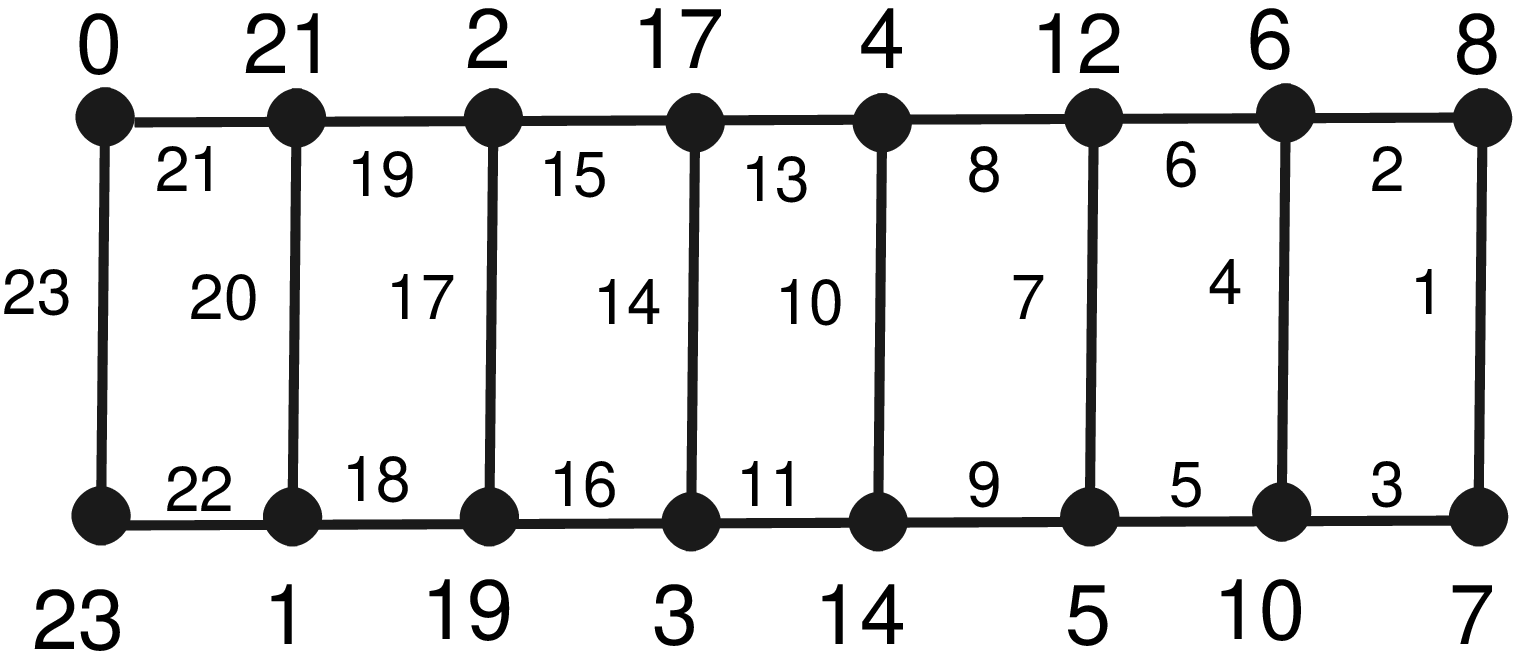}
\caption{$L_{16}$ }
\label{L16}
\end{center}
\end{figure}

\medskip
\noindent
Case 2: $k\equiv 3$(mod 4).\\
In this case $k=2t+1$ with $t\geq 5$ odd.
We define
$f: V(C_{2k})\rightarrow \{0,1,\ldots,2k+1\}$ in the following way
$$f(x_{2i+1})= \left\{
\begin{array}{lll}
i &\quad& \textrm{for}\ i \in \left[ 0,\frac{t-1}{2} \right]\\[2pt]
\frac{7t+5}{2}-i & \quad & \textrm{for}\ i \in \left[ \frac{t+1}{2},t\right]\\[2pt]
3t+2-i & \quad & \textrm{for}\ i \in \left[ t+1,t+\lfloor\frac{t+1}{4}\rfloor\right]\\[2pt]
t+2+i & \quad & \textrm{for}\ i \in \left[ t+\lfloor\frac{t+1}{4}\rfloor+1,\frac{3t-1}{2}\right]\\[2pt]
i-t & \quad & \textrm{for}\ i \in \left[ \frac{3t+1}{2},2t\right]
\end{array}
\right.$$
$$f(x_{2i})= \left\{
\begin{array}{lll}
4t+4-i &\quad& \textrm{for}\ i \in \left[1,\frac{t+1}{2} \right]\\[2pt]
\frac{t-1}{2}+i & \quad & \textrm{for}\ i \in \left[ \frac{t+3}{2},t+1\right]\\[2pt]
t+1+i & \quad & \textrm{for}\ i \in \left[ t+2,t+\lfloor\frac{t+3}{4}\rfloor\right]\\[2pt]
3t+2-i & \quad & \textrm{for}\ i \in \left[ t+\lfloor\frac{t+3}{4}\rfloor+1, \frac{3t+1}{2}  \right]\\[2pt]
5t+4-i & \quad & \textrm{for}\ i \in \left[ \frac{3t+3}{2},2t+1\right].
\end{array}
\right.$$
Also in this case one can directly check that
$$f({\cal O})=\left[0,\frac{t-1}{2}\right]\cup \left[\frac{5t+5}{2},3t+2\right] \cup
\left[2t+2-\left\lfloor\frac{t+1}{4}\right\rfloor, 2t+1\right]\cup$$
$$\cup \left[ 2t+2+\left\lfloor \frac{t+1}{4}\right\rfloor,\frac{5t+3}{2}  \right]
\cup \left[ \frac{t+1}{2},t \right]$$

$$f({\cal E})=\left[\frac{7t+7}{2},4t+3\right] \cup \left[t+1,\frac{3t+1}{2}\right] \cup
\left[2t+3,2t+1+\left\lfloor \frac{t+3}{4}\right\rfloor \right]\cup$$
$$\cup\left[\frac{3t+3}{2},2t+2-\left\lfloor\frac{t+3}{4}\right\rfloor\right]\cup \left[ 3t+3,\frac{7t+5}{2}\right].$$
By long calculations one can verify that, also in case 2, $f$ is a $2$-labeling of $C_{2k}$.
\hfill $\Box$\\
\\
In Figure \ref{C18-22} we have the $2$-graceful labeling of $C_{26}$ and that of $C_{30}$
obtained through the construction given in the proof of Theorem \ref{C2k}.
\\
\\
\noindent
To conclude we consider $d$-graceful $\alpha$-labelings for ladders.
We recall that the \emph{ladder graph} of order $2k$, denoted by $L_{2k}$,
can be seen as the cartesian product of the path $P_2=(0 \sim 1)$ of length $1$
by the path $P_k=(0 \sim 1 \sim \ldots \sim k-1)$ of length $k-1$. The ladder graph $L_{2k}$
is clearly bipartite with bipartite sets $A_0$, $A_1$ where $A_h$ is the set of pairs
$(i,j)\in \{0,1\}\times \{0,1,\ldots,k-1\}$ such that $i+j$ has the same parity as $h$.
\begin{thm}\label{ladder2}
The ladder graph $L_{2k}$ has a $2$-graceful $\alpha$-labeling if and only if $k$ is even.
\end{thm}
Proof.
Since $L_{2k}$ has $e=3k-2$ edges, $2$ divides $e$ if and only if $k$ is even.
Suppose now $k$ to be even. One can check that a $2$-graceful $\alpha$-labeling of $L_{2k}$ is given by
$$ f(i,j)=\left\{
\begin{array}{l}
j\ {\rm if} \ (i,j)\in A_0\\
3k-2j-1 \ {\rm if} \ (i,j)\in A_1,\ j=0,\ldots, \frac{k}{2}-1\\
3k-2j-2 \ {\rm if} \ (i,j)\in A_1,\ j=\frac{k}{2},\ldots, k-1.
\end{array}
\right.
$$
\hfill $\Box$\\
\\
We observe that the $2$-graceful $\alpha$-labeling constructed in the previous theorem
is very similar to the $\alpha$-labeling of $L_{2k}$ proposed by the author in \cite{AP}.\\
\\
Figure \ref{L16} shows the $2$-graceful $\alpha$-labeling
of $L_{16}$ provided by Theorem \ref{ladder2}.
\\
\\
As an immediate consequence of Theorem \ref{prop:alpha} we have
\begin{prop}
Let $k$ be an even integer. There exists a cyclic $L_{2k}$-decomposition of
$K_{\frac{3k}{2}\times 4n}$
for any integer $n\geq 1$.
\end{prop}

\noindent
\textbf{Acknowledgement.} The author is grateful to Marco Buratti
who suggested the item of this paper and to the anonymous referee
for his helpful and constructive comments.

\end{document}